\documentclass[a4paper, BCOR=0.0mm, DIV=calc]{scrartcl}
\usepackage[T1]{fontenc}
\usepackage[utf8]{inputenc}
\usepackage[ngerman]{babel}
\usepackage[osf,sc]{mathpazo}
\usepackage{textcomp}
\usepackage{microtype}
\usepackage{amsmath, amssymb, amsfonts}
\usepackage{tabularx}
\KOMAoptions{twoside=false, twocolumn=false, headinclude=false, footinclude=false, mpinclude=false, pagesize=auto}
\recalctypearea

\newcommand{\etc}{\text{etc}} 
\setkomafont{section}{\mdseries\scshape\Large}
\setkomafont{title}{\mdseries\scshape}

\begin{document}
\title{Über divergente Reihen}
\author{Leonhard Euler}
\date{}
\maketitle
\paragraph{§1}
Weil konvergente Reihen so definiert werden, dass sie aus stetig schrumpfenden Termen bestehen, die schließlich, wenn die Reihe ins Unendliche fortgeschritten ist, völlig verschwinden, sieht man leicht ein, dass die infinitesimalen Terme welcher Reihen nicht $0$ werden, sondern entweder endlich bleiben oder sogar ins Unendliche wachsen, dass diese Reihen dann, weil sie nicht konvergent sind, zur Klasse der divergenten Reihen gezählt werden müssen. Je nachdem, ob also die letzten Terme einer Reihe, zu welchen man, nachdem die Progression ins Unendliche fortgesetzt worden ist, von endlicher oder von unendlicher Größe waren, wird man zwei Arten der divergenten Reihen haben, von welchen jede der beiden weiter in zwei Gattungen unterteilt wird, je nachdem ob alle Terme mit denselben Vorzeichen versehen sind oder sich die Vorzeichen $+$ und $-$ abwechseln. Insgesamt werden wir also vier Gattungen der divergenten Reihen haben, aus welchen ich zur besseren Einsicht einige Beispiele hinzufügen möchte:
\[
\begin{array}{rl}
	\text{I.} & 1 + 1 + 1 + 1 + 1 + 1 + \etc \\[1ex]
	&\frac{1}{2} + \frac{2}{3} + \frac{3}{4} + \frac{4}{5} + \frac{5}{6} + \frac{6}{7} + \etc \\[1ex]
	\text{II.} & 1 - 1 + 1 - 1 + 1 - 1 + 1 - \etc \\[1ex]
	&\frac{1}{2} - \frac{2}{3} + \frac{3}{4} - \frac{4}{5} + \frac{5}{6} - \frac{6}{7} + \etc \\[1ex]
	\text{III.} & 1 + 2 + 3 + 4 + 5 + 6 + \etc \\[1ex]
	&1 + 2 + 4 + 8 + 16 + 32 + \etc \\[1ex]
	\text{IV.} & 1 - 2 + 3 - 4 + 5 - 6  + \etc \\[1ex]
	&1 - 2 + 4 - 8 + 16 - 32 + \etc 
\end{array}
\]
\paragraph{§2}
Über die Summen divergenter Reihen dieser Art herrscht große Uneinigkeit unter den Mathematikern, während die einen es verneinen, bestätigen die anderen, dass sie in einer Summe erfasst werden können. Und zuerst ist freilich klar, dass die Summen der Reihen, welche ich zur ersten Gattung gezählt haben, in der Tat unendlich groß sind, weil man, indem man die Terme tatsächlich zusammenfasst, zu einer Summe größer als jede gegebene Zahl gelangt; daher besteht natürlich kein Zweifel, dass die Summen dieser Reihen durch Ausdrücke dieser Art $\frac{a}{0}$ beschafft werden können. Die Diskussion unter den Geometern geht also hauptsächlich über die übrigen Gattungen, und die Argumente, die auf beiden Seiten zur Verteidigung der eigenen Äußerungen angeführt werden, haben eine so große Überzeugungskraft, dass noch keine Partei gezwungen werden konnte, der anderen Recht zu geben.
\paragraph{§3}
Aus der zweiten Gattung hat Leibniz als erster diese Reihe betrachtet
\[
	1 - 1 + 1 - 1 + 1 - 1 + 1 - 1  + \etc,
\]
deren Summe er festgesetzt hatte gleich $\frac{1}{2}$ zu sein, und hatte das auf diese hinreichend strenge Begründung gestützt: Zuerst geht diese Reihe hervor, wenn dieser Bruch $\frac{1}{1+a}$ durch wiederholte Teilung auf gewohnte Weise in diese Reihe aufgelöst wird
\[
	1 - a + a^2 - a^3 + a^4 - a^5 + \etc
\] 
und der Wert des Buchstaben $a$ der Einheit gleich genommen wird. Um dann aber dies noch mehr zu bestätigen und die, die an die Rechnung nicht gewöhnt sind, zu überzeugen, gebrauchte er die folgende Begründung: Wenn die Reihe irgendwo abgebrochen wird und die Anzahl der Terme eine gerade Zahl war, dann wird ihr Wert gleich $0$ sein, wenn aber die Anzahl der Terme ungerade war, wird der Wert der Reihe gleich $1$ sein; wenn die Reihe daher also ins Unendliche fortschreitet und die Anzahl der Terme weder als gerade noch als ungerade angesehen werden kann, schloss er, dass die Summe weder gleich $0$ noch gleich $1$ sein kann, sondern einen bestimmten Mittelwert von jenen beiden gleich verschieden haben muss, welcher gleich $\frac{1}{2}$ ist.
\paragraph{§4}
Gegen diese Argumente pflegt von den Gegnern folgendes entgegnet zu werden: Zuerst sei der Bruch $\frac{1}{1+a}$ nicht der unendlichen Reihe
\[
	1 - a + a^2 - a^3 + a^4 - a^5 + a^6 - \etc
\]
gleich, wenn $a$ nicht ein Bruch kleiner als die Einheit ist. Wenn nämlich die Teilung irgendwann abgebrochen wird und dem Quotienten vom Rest der entsprechende Teil hinzugefügt wird, werde es offensichtlich zu Widersprüchen führen; denn es werde
\[
	\frac{1}{1+a} = 1 - a + a^2 - a^3 + a^4 - \dots \pm a^n \mp \frac{a^{n+1}}{1+a},
\]
sein, und obwohl die Zahl $n$ als unendlich festgesetzt wird, lässt sich der Bruch $\mp \frac{a^{n+1}}{1+a}$ nicht weglassen, außer wenn er tatsächlich verschwindet, was nur in den Fällen passiert, in denen $a<1$ ist, und die Reihe konvergent wird. In den übrigen Fällen müsse aber diese Mantisse $\mp \frac{a^{n+1}}{1+a}$ immer beibehalten werden, und obwohl sie mit dem zweifelhaften Vorzeichen $\pm$, je nachdem ob $n$ eine gerade Zahl oder eine ungerade Zahl war, versehen ist, könne sie dennoch, wenn $n$ unendlich ist, daher nicht weggelassen werden, weil eine unendliche Zahl weder gerade noch ungerade sei und man deshalb keinen Grund habe, ob eines der beiden Vorzeichen eher zu verwenden ist; es sei nämlich absurd zu glauben, dass eine ganze Zahl gegebenen ist, nicht einmal eine unendliche, die weder gerade noch ungerade sein soll.
\paragraph{§5}
Aber bei dieser Entgegnung pflegt von jenen, die den divergenten Reihen bestimmte Summen zuteilen, mit Recht getadelt zu werden, dass die unendliche Zahl als eine bestimmte Zahl verstanden wird und sogar als entweder gerade oder ungerade festgesetzt wird, obwohl sie dennoch unbestimmt ist. Denn sofort werde auch eine Reihe bezeichnet ins Unendliche fortzuschreiten, entgegen der eigentlichen Idee, wenn ein Term derselben Reihe als letzter, wenn auch infinitesimaler, verstanden wird; und daher löse sich die zuvor erwähnte Entgegnung über die dem letzten Term hinzuzufügende oder abzunehmende Mantisse von selbst auf. Weil man also in einer unendlichen Reihe niemals zum Ende gelangt, könne man auch niemals zu einer Stelle solcher Art gelangen, wo es notwendig wäre, jene Mantisse hinzuzufügen, und diese könne daher nicht nur missachtet werden, sondern müsse es sogar, weil ihr niemals ein Platz überlassen wird. Und diese Argumente, die, um die Summen divergenter Reihen zu bestätigen oder zu widerlegen, angeführt werden, beziehen sich auch immer auf die vierte Gattung, die außerdem durch keine neuen Argumente hin und her zu diskutiert werden pflegt.
\paragraph{§6}
Aber die, die gegen die Summen der divergenten Reihen plädieren, glauben in der dritten Gattung den sichersten Schutz zu finden. Denn obwohl die Terme dieser Reihen immer weiter wachsen und man daher durch tatsächliches Zusammenfassen der Terme zu Summen größer als jede angebbare Zahl gelangen kann, was die Definition des Unendlichen ist, sind die Verteidiger der Summen gezwungen in dieser Gattung Reihen zuzulassen, deren Summen endlich und sogar konvergent oder kleiner als $0$ sind. Weil nämlich der Bruch $\frac{1}{1-a}$ durch Teilung in eine Reihe entwickelt
\[
	1 + a + a^2 + a^3 + a^4 + \etc
\]
gibt, müsste
\begin{align*}
-1 &= 1 + 2 + 4 + 8 + 16 + \etc \\
-\frac{1}{2} &= 1 + 3 + 9 + 27 + 81 + \etc
\end{align*}
sein, was den Gegnern nicht zu Unrecht als höchstgradig absurd erscheint, weil man durch Addition positiver Zahlen niemals zu negativen Summen gelangen kann. Und daher drängen sie umso mehr auf die Notwendigkeit der hinzuzufügenden Mantisse, die zuvor erwähnt wurde, weil, nachdem diese hinzugefügt wurde, klar ist, dass
\[
	-1 = 1 + 2 + 4 + 8 + \dots + 2^n + \frac{2^{n+1}}{1-2}
\]
ist, auch wenn $n$ eine Unendliche Zahl ist.
\paragraph{§7}
Die Fürsprecher der Summen divergenter Reihen setzen daher, um dieses außergewöhnliche Paradoxon zu erklären, einen feinen, mehr als wahren Unterschied zwischen den negativen Größen fest; sie erklären, dass, während die einen kleiner als $0$ sind, die anderen aber größer als unendlich oder mehr als unendlich sind. Denn einen Wert "`$-1$"' müsse man natürlich nehmen, wann immer er verstanden wird aus der Subtraktion einer größeren Zahl $a+1$ von einer kleineren $a$ zu entspringen, den anderen aber, wann immer er jener Reihe $1+2+4+8+16+\etc$ gleich gefunden wird und aus der Teilung der Zahl $+1$ durch $-1$ entsteht; in jenem Fall sei ja die Zahl kleiner als $0$, in diesem aber größer als unendlich. Für mehr Bestätigung führen sie dieses Beispiel der Reihe der Brüche an
\[
	\frac{1}{4},\quad \frac{1}{3},\quad \frac{1}{2},\quad \frac{1}{1},\quad \frac{1}{-1},\quad \frac{1}{-2},\quad \frac{1}{-3},\quad \etc,
\]
welche, weil sie bei den ersten Termen als wachsend erkannt wird, auch notwendigerweise anzusehen ist, immer zu wachsen, woher sie folgern, dass $\frac{1}{-1} > \frac{1}{0}$ und $\frac{1}{-2} > \frac{1}{-1}$ ist und so weiter; und daher, insoweit $\frac{1}{-1}$ durch $-1$ und $\frac{1}{0}$ durch unendliche "`$\infty$"' ausgedrückt wird, dass $-1 > \infty$ ist; auf diese Weise weisen sie jene sich ergebende Absurdität genügt geistreich von sich.
\paragraph{§8}
Obwohl aber diese Unterscheidung geistreich erdacht scheint, genügt sie den Gegnern wenig, und scheint sogar sichere Regeln der Analysis zu verletzen. Denn wenn nämlich jene zwei Werte von $-1$, insoweit sie entweder gleich $1-2$ oder gleich $\frac{1}{-1}$ sind, sich in der Tat voneinander unterscheiden, dass sich diese nicht zusammenbringen lassen, würde die Sicherheit und der Gebrauch der Regeln, denen wir beim Kalkül folgen, völlig aufgehoben werden, was gewiss noch absurder wäre als das, weswegen diese Unterscheidung erdacht worden ist; wenn aber $1-2 = \frac{1}{-1}$ ist, wie es die Vorschriften der Algebra erfordern, wird die Aufgabe keineswegs erledigt, weil jene Größe $-1$, die der Reihe $1+2+4+8+\etc$ gleich gesetzt wird, kleiner als nichts ist und daher dieselbe Schwierigkeit zurückbleibt. Dennoch scheint es wiederum mit der Wahrheit vereinbar, wenn wir sagen, dass dieselben Größen, die kleiner sein sollen als nichts, zugleich größer als unendlich werden können. Denn nicht nur aus der Algebra, sondern auch aus der Geometrie lernen wir, dass zwei Übergänge von positiven zu negativen Größen gegeben sind, zum einen durch $0$ oder Nichts, zum anderen durch das Unendliche, und sogar, dass die Größen, indem sie von der $0$ wachsen und schrumpfen, aufeinander zurückgehen und schließlich zum selben Term $0$ zurückkehren, sodass die Größen größer als die kleineren als $0$ und die Größen kleiner als unendlich mit den Größen größer als Nichts übereinstimmen.
\paragraph{§9}
Dieselben aber verneinten, dass diese Summen der divergenten Reihen, die angegeben zu werden pflegen, richtig sind, bringen nicht nur keine anderen hervor, sondern beschließen auch, sich darum zu bemühen, dass die Summe der divergenten Reihe nur ausgedacht ist. Denn sie könnten die Summe der konvergenten Reihen wie z.\,B. dieser
\[
	1 + \frac{1}{2} + \frac{1}{4} + \frac{1}{8} + \frac{1}{32} + \etc
\]
nur daher als gleich $2$ zulassen, weil, je mehr Terme dieser Reihe wir tatsächlich addieren, wir umso näher an die $2$ gelangen; bei divergenten Reihen verhalte sich aber die Sache bei weitem nicht so; je mehr Terme wir nämlich addieren, umso mehr unterscheiden sich die Summen, die hier hervorgehen, voneinander und gehen nie zu einem festen und bestimmten Wert heran. Daher folgern sie, dass nicht einmal die Idee einer Summe auf divergente Reihen übertragen werden kann und die Mühe derer, die beim Untersuchen der Summen divergenter Reihen aufgebracht wird, natürlich unnötig ist und den wahren Prinzipien der Analysis entgegen ist.
\paragraph{§10}
So real diese Uneinigkeit aber auch erscheinen mag, kann dennoch keine Partei von der anderen eines Fehlers überführt werden, sooft in der Analysis der Gebrauch von Reihen dieser Art auftaucht; es muss von großer Bedeutung sein, dass keine Partei falsch liegt, sondern die Uneinigkeit in den Werten allein gelegen ist. Wenn ich nämlich bei einer Rechnung zu dieser Reihe $1-1+1-1+1-1+\etc$ gelange und an ihrer Stelle $\frac{1}{2}$ einsetze, wird gewiss niemand mir mit Recht einen Fehler anlasten, der dennoch jedem sofort ins Auge spränge, wenn ich irgendeine andere Zahl an deren Stelle gesetzt hätte; daher kann kein Zweifel bestehen bleiben, dass die Reihe $1-1+1-1+1-1+\etc$ und der Bruch $\frac{1}{2}$ äquivalente Größen sind, und sich die eine anstelle der anderen immer ohne einen Fehler einsetzen lässt. Die ganze Frage scheint also scheint nur darauf zurückzugehen, ob wir den Bruch $\frac{1}{2}$ richtigerweise die Summe der Reihe $1-1+1-1+\etc$ nennen; die das hartnäckig verneinen, obwohl sie dennoch nicht wagen die Äquivalenz zu verneinen, sind dafür vehement zu verachten, nicht die Logik zu beachten.
\paragraph{§11}
Ich glaube aber, dass dieser Streit beigelegt werden wird, wenn wir uns eifrig auf das folgende beziehen wollten. Sooft wir in der Analysis zu einem gebrochenen oder transzendenten Ausdruck gelangen, pflegen wir diesen genauso oft in eine geeignete Reihe zu verwandeln, auf welche die folgende Rechnung gefälliger angewendet werden kann. So weit also nur unendliche Reihen in der Analysis Platz finden, so weit sind die aus der Entwicklung eines endlichen Ausdrucks entsprungen, und deswegen lässt sich in der Rechnung immer anstelle der unendlichen Reihe die Form einsetzen, aus deren Entwicklung sie entstanden ist; Wie daher mit dem größten Ertrag die Regeln angegeben zu werden pflegen, endliche Ausdrücke, die aber mit einer weniger geeigneteren Form versehen sind, in unendliche Reihen zu verwandeln, so sind umgekehrt die Regeln als die nützlichsten anzusehen, mit deren Hilfe, wenn irgendeine unendliche Reihe vorgelegt war, der endliche Ausdruck gefunden werden kann, aus welchem sie resultiert. Und weil dieser Ausdruck immer ohne Fehler anstelle der unendlichen Reihe eingesetzt werden kann, ist es notwendig, dass der Wert jeder der beiden derselbe ist; daraus erreicht man, dass keine unendliche Reihe gegeben ist, dass nicht auch gleichzeitig der endliche Ausdruck als jener äquivalent angesehen werden kann.
\paragraph{§12}
Wenn wir also die gewohnte Bezeichnung der Summe nur so ändern, dass wir sagen, dass die Summe der endliche Ausdruck der Reihe ist, aus dessen Entwicklung die Reihe selbst entsteht, werden alle Schwierigkeiten, die von beiden Parteien erwähnt worden sind, von selbst verschwinden. Zuerst beschafft nämlich der Ausdruck, aus dessen Entwicklung die konvergente Reihe entspringt, zugleich seine Summe, wobei sie diese Bezeichnung dann im gewöhnlichen Sinne erhalten hat, und, wenn die Reihe divergent war, kann die Frage nicht weiter als absurd bezeichnet werden, wenn wir sie als endlichen Ausdruck untersuchen, der nach den analytischen Regeln entwickelt jene Reihe selbst erzeugt. Und weil sich dieser Ausdruck in der Rechnung anstelle der Reihe einsetzen lässt, werden wir nicht bezweifeln können, dass er derselben gleich ist. Nachdem das erklärt worden ist, wollen wir nicht einmal vom gewohnten Sprachgebrauch abweichen, wenn wir den Ausdruck, der der Reihe gleich ist, auch als Summe bezeichnen, solange wir für die divergenten Reihen die Benennung nicht mit der Idee der Summe verbinden, weil, je mehr Terme tatsächlich gesammelt werden, man umso näher an den wahren Wert der Summe herangehen müsste.
\paragraph{§13}
Nachdem diese Dinge vorausgeschickt worden sind, glaube ich, dass es keinen geben wird, der mich als zu Tadelnden ansieht, weil ich im folgenden gründlicher die Summe der Reihe
\[
	1 - 1 + 2 - 6 + 24 - 120 + 720 - 5040 + 40320 - \mathrm{etc}
\]
untersuchen werde, welche Reihe von Wallis hypergeometrisch genannt wurde, hier ist sie mit alternierenden Vorzeichen versehen worden. Diese Reihe scheint aber umso bemerkenswerter, weil ich hier viele Summierungsmethoden, die mir anderenorts bei einer Aufgabe dieser Art einen riesigen Nutzen geleistet haben, vergeblich ausprobiert habe. Zuerst lässt sich freilich zweifeln, ob diese Reihe eine endlich Summe hat oder nicht, weil sie um Vieles mehr divergiert als eine einzige geometrische Reihe; dass die Summe der geometrischen Reihen aber endlich ist, ist außer Zweifel gestellt worden; aber weil doch bei den geometrischen die Divergenz nicht dagegen spricht, dass sie summierbar sind, so scheint es wahrscheinlich, dass auch diese hypergeometrische Reihe eine endliche Summe hat. Es wird also in Zahlen, zumindest näherungsweise, der Wert ihres endlichen Ausdruckes gesucht, aus dessen Entwicklung die vorgelegte Reihe selbst entsteht.
\paragraph{§14}
Zuerst aber habe ich die Methode benutzt, die auf diesen Fundament beruht: wenn eine Reihe dieser Art vorgelegt ist
\[
	s = a - b + c - d + e - f + g - h + \mathrm{etc}
\]
und, nachdem die Vorzeichen der Terme $a$, $b$, $c$, $d$, $e$, $f$, etc weggelassen worden sind, die Differezen
\[
	b-a,\quad c-b,\quad d-c,\quad e-d,\quad \mathrm{etc}
\]
genommen werden und weiter deren Differenzen
\[
	c-2b+a,\quad d-2c+b,\quad e-2d+c,\quad \mathrm{etc},
\]
die zweite Differenzen genannt werden, und nach dem gleichen Bildungsgesetz die dritten, vierten, fünften Differenzen etc genommen werden, dann sage ich, wenn die Terme dieser ersten, zweiten, dritten, vierten Differenzen $\alpha$, $\beta$, $\gamma$, $\delta$ etc sind, dass die Summe der vorgelegten Reihe selbst
\[
	s = \frac{a}{2} - \frac{\alpha}{2} + \frac{\beta}{8} - \frac{\gamma}{16} + \frac{\delta}{32} - \mathrm{etc}
\]
sein wird, welche, wenn sie nicht schon konvergent ist, dennoch gewiss um Vieles mehr konvergieren wird als die vorgelegte; daher wird man, wenn dieselbe Methode erneut auf diese letztere Reihe angewandt wird, den Wert oder die Summe $s$ durch eine noch stärker konvergente Reihe ausgedrückt finden.
\paragraph{§15}
Diese Methode hat den größten Nutzen beim Summieren divergenter Reihen der zweiten und vierten Art, ob man nun schließlich zu konstanten Differenzen gelangt oder ob anders, sofern die Divergenz nicht allzu groß ist. Wenn so
\[
	s = 1 - 1 + 1 - 1 + 1 - \text{etc}
\]
ist, wird wegen
\[
	a = 1,\quad \alpha = 0,\quad \beta = 0,\quad \text{etc}
\]
gleich
\[
	s = \frac{1}{2}
\]
sein. Wenn
\[
\begin{array}{rlcccccccccccccc}
	s &=& 1 &-& 2 &+& 3 &-& 4 &+& 5 &-& 6 &+& \text{etc} \\
	\text{diff\,I.} & & &1& &1& &1& &1& &1&
\end{array}
\]
ist, wird
\[
	s = \frac{1}{2} - \frac{1}{4} = \frac{1}{4}
\]
sein, wie schon anderswoher hinreichend bekannt ist. Wenn
\[
\begin{array}{rlcccccccccccccc}
	s &=& 1 &-& 4 &+& 9 &-& 16 &+& 25 &-& 36 &+& \text{etc} \\
	\text{diff\,I.} & & &3& &5& &7& &9& &11& \\
	\text{diff\,II.} & & & &2& &2& &2& &2& &
\end{array}	
\]
ist, wird
\[
	s = \frac{1}{2} - \frac{3}{4} + \frac{2}{8} = 0
\]
sein, wie auch bekannt ist. Wenn
\[
\begin{array}{rlcccccccccccccc}
	s &=& 1 &-& 3 &+& 9 &-& 27 &+& 81 &-& 243 &+& \text{etc} \\
	\text{diff\,I.}   & & &2& &6& &18& &54& &162& \\
	\text{diff\,II.}  & & & &4& &12& &36& &108& & \\
	\text{diff\,III.} & & & & &8& &24& &72& & & \\
	\text{diff\,IV.}  & & & & & &16& &48& & & & \\
	                  & & & & & &  & \text{etc}&  & & & & 
\end{array}	
\]
ist, wird
\[
	s = \frac{1}{2} - \frac{2}{4} + \frac{4}{8} - \frac{8}{16} + \text{etc} = \frac{1}{2} - \frac{1}{2} + \frac{1}{2} - \frac{1}{2} + \text{etc} = \frac{1}{4}
\]
sein.
\paragraph{§16}
Man wende gleich diese Methode auf die vorgelegte Reihe an
\[
	A = 1 - 1 + 2 - 6 + 24 - 120 + 720 - 5040 + 40320 - \etc,
\]
die wegen $1-1 = 0$, wenn sie durch $2$ geteilt wird, übergeht in diese
\begin{center}
$\frac{A}{2} = 1 - 3 + 12 - 60 + 360 - 2520 + 20160 - 181440 + \etc$\\
$2,\quad 9,\quad 48,\quad 300,\quad 2160,\quad 17640,\quad 161280$\\
$7,\quad 39,\quad 252,\quad 1860,\quad 15480,\quad 143640$ \\
$32,\quad 213,\quad 1608,\quad 13620,\quad 128160$\\
$181,\quad 1395,\quad 12012,\quad 114540$\\
$1214,\quad 10617,\quad 102528$\\
$9403,\quad 91911$\\
$8250$\\
\end{center}
Daher also folgt, dass
\[
	\frac{A}{2} = \frac{1}{2} - \frac{2}{4} + \frac{7}{8} - \frac{32}{16} + \frac{181}{32} - \frac{1214}{64} + \frac{9403}{128} - \frac{82508}{256} + \etc
\]
sein wird oder
\begin{center}
	$A = \dfrac{7}{4} - \dfrac{32}{8} + \dfrac{181}{16} - \dfrac{1214}{32} + \dfrac{9403}{64} - \dfrac{82508}{128} + \etc$\\[1ex]
	$\dfrac{18}{8},\quad \dfrac{117}{16},\quad \dfrac{852}{32},\quad \dfrac{6975}{64},\quad \dfrac{63702}{128}$\\[1ex]
	$\dfrac{81}{16},\quad \dfrac{618}{32},\quad \dfrac{5271}{64},\quad \dfrac{49752}{128}$\\[1ex]
	$\dfrac{456}{32},\quad \dfrac{4035}{64},\quad \dfrac{39210}{128}$\\[1ex]
	$\dfrac{3123}{64},\quad \dfrac{31140}{128}$\\[1ex]
	$\dfrac{24894}{128}$
\end{center}
Also
\[
	A = \frac{7}{8} - \frac{18}{32} + \frac{81}{128} - \frac{456}{512} + \frac{3123}{2048} - \frac{24894}{8192} + \etc
\]
oder
\begin{center}
	$A - \dfrac{5}{16} = \dfrac{81}{128} - \dfrac{456}{512} + \dfrac{3123}{2048} - \dfrac{24894}{8192} + \etc$ \\[1ex]
	$\quad\dfrac{132}{512},\quad \dfrac{1299}{2048},\quad \dfrac{12402}{8192}$\\[1ex]
	$\quad\dfrac{771}{2048},\quad \dfrac{7206}{8192}$\\[1ex]
	$\quad\dfrac{4122}{8192}$
\end{center}
Also
\[
	A - \frac{5}{16} = \frac{81}{256} - \frac{132}{2048} + \frac{771}{16384} - \frac{4122}{131072}
\]
oder
\[
	A = \frac{5}{16} + \frac{516}{2048} + \frac{2046}{131072} + \etc = \frac{38015}{65536} = 0,580.
\]
Es scheint also die Summe dieser Reihe fast gleich $0,580$ zu sein; wegen der weggelassenen Terme aber wird sie ein wenig größer sein, was überaus mit dem unten zu zeigenden übereinstimmt, wo die Summe dieser Reihe gezeigt werden wird gleich $0,5963736$ zu sein; zugleich aber ist klar, dass diese Methode nicht hinreichend geeignet ist um die Summe so exakt zu bestimmen.
\paragraph{§17}
Darauf bin ich auf andere Weise die Sache so angegangen: Es sei diese Reihe vorgelegt
\[
\begin{array}{ccccccccccc}
	& 1 & 2 & 3 & 4 & 5 & 6 & 7 & \dots & n & n+1 \\
	B  & 1, & 2, & 5, & 16, & 65, & 326, & 1957, & \dots & P, & nP+1 
\end{array}
\]
die Differenzen sind
\begin{center}
$1,\quad 3,\quad 11,\quad 49,\quad 261,\quad 1631$ \\
$2,\quad 8,\quad 38,\quad 212,\quad 1370$ \\
$6,\quad 30,\quad 174,\quad 1158$ \\
$24,\quad 144,\quad 984$ \\
$120,\quad 840$ \\
$720$
\end{center}
weil die ersten Terme dieser weiteren Differenzen davon
\[
	1,\quad 2,\quad 6,\quad 24,\quad 120,\quad 720,\quad \etc
\]
sind, wird der dem Index $n$ entsprechende Term
\begin{align*}
	P =& 1 + (n-1) + (n-1)(n-2) + (n-1)(n-2)(n-3) \\
	&+ (n-1)(n-2)(n-3)(n-4) + \etc
\end{align*}
sein. Daher wird, wenn $n=0$ wird, der dem Index $n$ entsprechende Term oder der als erste vorangehende gleich
\[
	1 - 1 + 2 - 6 + 24 - 120 - \etc = A
\]
sein, sodass, wenn der dem Index $n$ entsprechende Term dieser Reihe gefunden werden könnte, derselbe zugleich der Wert oder die Summe der vorgelegten Reihe
\[
	A = 1 - 1 + 2 - 6 + 24 - 120 + 720 - \etc
\]
sein würde. Wenn daher also jene Reihe $B$ invertiert wird, dass man die Reihe
\[
\begin{array}{ccccccccccc}
	& 1 & 2 & 3 & 4 & 5 & 6 & 7 \\[1ex]
	C & 1, & \frac{1}{2}, & \frac{1}{5}, & \frac{1}{16}, & \frac{1}{65}, & \frac{1}{326}, & \frac{1}{1957}, & \etc 
\end{array}
\] 
hat, wird der dem Index $n$ entsprechende Term dieser Reihe gleich $\frac{1}{A}$ sein, woher aus ihm auch der Wert von $A$ selbst erkannt werden können wird. Die einzelnen Differenzen dieser Reihe mögen mit den Termen $\alpha$, $\beta$, $\gamma$, $\delta$, $\varepsilon$, $\etc$ beginnen, mit natürlich hier so zu nehmenden Differenzen, dass je der Term vom vorhergehenden abgezogen wird; es wird der dem Index $n$ entsprechende Term
\[
	\frac{1}{P} = 1 - (n-1)\alpha + \frac{(n-1)(n-2)}{1\cdot 2}\beta - \frac{(n-1)(n-2)(n-3)}{1\cdot 2\cdot 3}\gamma + \etc
\]
sein. Und daher wird für $n=0$ gesetzt durch eine gewisse konvergente Reihe
\[
	\frac{1}{A} = 1 + \alpha + \beta + \gamma + \delta + \etc
\]
sein. Es ist aber, indem man diese Brüche in Dezimale entwickelt:
\[
\begin{array}{rccccc}
& \text{Diff. 1} & \text{Diff. 2} & \text{Diff. 3} & \text{Diff. 4} & \text{Diff. 5} \\
1 = 1,0000000 & & & & & \\
& 5000000 & & & & \\
\frac{1}{2} = 0,5000000 &  & 2000000 & & & \\ 
& 3000000 & & 375000 & & \\
\frac{1}{5} = 0,2000000 & & 1625000 & & -346154 & \\
& 1375000 & & 721154 & & -511445 \\
\frac{1}{16} = 0,0625000 & & 903848 & & +165291 & \\
& 471154 & & 555863 & & -140195\\
\frac{1}{65} = 0,0153846 & & 347983 & & +305486 & \\
& 123171 & & 250377 & & +131530\\
\frac{1}{326} = 0,0030675 & & 97606 & & +173956 & \\
& 25565 & & 76421 & & +114979\\
\frac{1}{1957} =0,0005110 & & 21185 & & +58977 & \\
& 4380 & & 17444 & & +44716\\
0,0000370 & & 3741 & & +14261 & \\
& 639 & & 3183 & & +11564\\
0,0000091 & & 558 & & +2697 & \\
& 81 & & 486 & & +2275\\
0,0000010 & & 72 & & +422 & \\
& 9 & & 64 & & +365\\
0,0000001 & & 8 & & +57 &  
\end{array}
\]
Aus diesen Differenzen also wird
\[
	\frac{1}{A} = 1,6517401 \quad \text{und} \quad A = 0,6
\]
sein, welcher Wert gut mit dem zuvor gefundenen zusammenpasst; aber dennoch ist diese Methode wegen der vierten Differenzen, der fünften und einigen der folgenden negativen nicht hinreichend sicher.
\paragraph{§18}
Wir wollen nun die Logarithmen der einzelnen Terme der Reihe $B$ nehmen, dass man diese neue Reihe hat:
\[
\begin{array}{ccccccccccl}
	& 1 & 2 & 3 & 4 & 5 & 6 & 7 & 8 \\
	D & \log{1}, &\log{2}, & \log{5}, & \log{16}, & \log{65}, & \log{326}, & \log{1957}, & \log{13700},& \etc 
\end{array}
\]
in welcher, nachdem die benachbarten Differenzen auf gewohnte Weise genommen worden sind, die ersten Terme $\alpha$, $\beta$, $\gamma$, $\delta$, $\varepsilon$, etc seien, und es wird der dem Index $0$ entsprechende Term dieser Reihe gleich
\[
	0 - \alpha + \beta - \gamma + \delta - \varepsilon + \etc
\]
sein, welcher also der Logarithmus der gesuchten Summe gleich $A$ ist. Es sind in der Tat diese Logarithmen mit den benachbarten Differenzen die folgenden:
\[
\begin{array}{ccccccccc}
& \text{Diff. 1} & \text{Diff. 2} & \text{Diff. 3} & \text{Diff. 4} & \text{Diff. 5} & \text{Diff. 6} & \text{Diff. 7} & \text{Diff. 8} \\
0,0000000 & & & & & & & & \\
& 0,3010300 & & & & & & & \\
0,3010300 & & 969100 & & & & & & \\
& 0,3979400 & & 103000 & & & & & \\
0,6989700 & & 1072100 & & -138666& & & & \\
& 0,5051500 & & -35666 & & +53006& & & \\
1,2041200 & & 1036434 & & -85660 & & +19562 & & \\
& 0,6087934 & & -121326 & & +72568 & & -57744 & \\
1,8129134 & & 915108 & & -12092 & & -38182 & & +65446 \\
& 0,7003042 & & -134418 & & +34386 & & +7702 & \\
2,5132176 & & 780690 & & +21294 & & -30480 & & \\
& 0,7783732 & & -113124 & & +3906 & & & \\
3,2915908 & & 667566 & & +25200 & & & & \\
& 0,8451298 & & -87925 & & & & & \\
4,1367206 & & 579641 & & & & & & \\
& 0,9030939 & & & & & & & \\
5,0398145 & & & & & & & &
\end{array}
\]
also wird
\[
\begin{array}{rcccccc}
& \text{Diff. 1} & \text{Diff. 2} & \text{Diff.3} & \text{Diff. 4} & \text{Diff. 5} & \text{Diff. 6} \\
\log{A} = & & & & & & \\  
-0,3010300 & & & & & &\\
& +2041200 & & & & & \\
+0,0969100 & & +1175100 & & & & \\
& +866100 & & +550666 & & & \\
-0,0103000 & & +624434 & & +359570 & & \\
& +241666 & & +191096 & & +826928 & \\
-0,0138666 & & +433338 & & -467358 & & +2133994 \\
& -191672 & & +658454 & & -1307066 & \\
-0,0053006 & & -225116 & & +839708 & & -2083670\\
& +33444 & & -181254 & & +776604 & \\
+0,0019562 & & -43862 & & +63103 & & \\
& +77306 & & -244357 & & & \\
+0,0057744 & & +200495 & & & & \\
& -123189 & & & & & \\
+0,0065445 & & & & & &
\end{array}
\]
sein, woher durch die zuvor erläuterte Methode
\[
	\log{\frac{1}{A}} = \frac{0,0310300}{2} + \frac{2041200}{4} + \frac{1175100}{8} + \frac{550666}{16} + \frac{359570}{32} + \frac{826928}{64} + \etc
\]
sein oder
\[
	\log{\frac{A}{1}} = 0,7779089 \quad \text{und daher} \quad A = 0,59966,
\]
welche Zahl sich aber noch berechnen lässt, leicht zu klein zu sein. Trotzdem kann auch auf diese Weise weder sicher noch hinreichend angenehm zur Erkenntnis des Wertes $A$ gelangt werden, auch wenn diese Methode unendliche Wege liefert, diesen Wert zu untersuchen; für diesen Zweck scheinen die einen derer um Vieles geeigneter als die anderen.
\paragraph{§19}
Wir wollen nun auch analytisch den Wert dieser Reihe untersuchen, wir wollen sie in der Tat in einem weiteren Sinne auffassen; es sei also
\[
	s = x - 1x^2 + 2x^3 - 6x^4 + 24x^5 - 120x^5 + \etc,
\]
welche differenziert
\[
	\frac{\mathrm{d}s}{\mathrm{d}x} = 1 - 2x + 6xx - 24x^3 + 120x^4 - \etc  = \frac{x-s}{xx}
\]
geben wird, woher
\[
	\mathrm{d}s + \frac{s\mathrm{d}x}{xx} = \frac{\mathrm{d}x}{x}
\]
wird, das Integral welcher Gleichung, wenn $e$ für die Zahl genommen wird, bei der der hyperbolischer Logarithmus gleich $1$ ist,
\[
	e^{-1:x}s = \int{\frac{e^{-1:x}}{x}\mathrm{d}x} \quad \text{und} \quad s = e^{1:x}\int{\frac{e^{-1:x}}{x}\mathrm{d}x}
\]
sein wird. In dem Fall also, in dem $x=1$ ist, wird
\[
	1 - 1 + 2 - 6 + 24 - 120 + \etc = e\int{\frac{e^{-1:x}}{x}\mathrm{d}x}
\]
sein. Es drückt also diese Reihe die Fläche unter der Kurve aus, deren Gestalt zwischen der Abszisse $x$ und $y$ in dieser Gleichung 
\[
	y = \frac{e\cdot e^{-1:x}}{x}
\]
enthalten ist, wenn die Abszisse gleich $1$ gesetzt wird, oder es wird
\[
	y = \frac{e}{e^{1:x}\cdot x}
\]
sein. Diese Kurve aber ist so beschaffen, dass für $x=0$ gesetzt $y=0$ wird; wenn aber $x=1$ ist, wird $y=1$ sein; die Mittelwerte der Ordinaten werden sich in der Tat so verhalten, dass
\[
\begin{array}{lcllcl}
\text{wenn galt} & \qquad \qquad & \text{dann auch wird} \qquad & \text{wenn galt} & \qquad \qquad & \text{dann auch wird} \\[3mm]
x = \dfrac{0}{10} & \qquad \qquad & y = 0                 \qquad & x = \dfrac{5}{10} & \qquad \qquad & y = \dfrac{10}{5e^{5:5}} = \dfrac{2}{e} \\[3mm]
x = \dfrac{1}{10} & \qquad \qquad & y = \dfrac{10}{e^{9:1}} \qquad & x = \dfrac{6}{10} & \qquad \qquad & y = \dfrac{10}{6e^{4:6}} \\[3mm]
x = \dfrac{2}{10} & \qquad \qquad & y = \dfrac{10}{2e^{8:2}} \qquad & x = \dfrac{7}{10} & \qquad \qquad & y = \dfrac{10}{7e^{3:7}} \\[3mm]
x = \dfrac{3}{10} & \qquad \qquad & y = \dfrac{10}{3e^{7:3}} \qquad & x = \dfrac{8}{10} & \qquad \qquad & y = \dfrac{10}{8e^{2:8}} \\[3mm]
x = \dfrac{4}{10} & \qquad \qquad & y = \dfrac{10}{4e^{6:4}} \qquad & x = \dfrac{9}{10} & \qquad \qquad & y = \dfrac{10}{9e^{1:9}}
\end{array}
\]
Nachdem also diese Kurve konstruiert wurde, wird sofort klar werden, dass deren Fläche, die der Abszisse $x=1$ entspricht, nicht nur endlich ist, sondern sogar kleiner ist als das Quadrat mit Seitenlängen gleich $1$, größer aber als dessen Hälfte $\frac{1}{2}$. Wenn also tatsächlich die Grundseite $x=1$ in $10$ gleiche Teile geteilt wird und die Anteile der Fläche als Trapeze betrachtet werden und man die Flächen untersucht, wird man den wahren Wert der Reihe
\[
	1 - 1 + 2 - 6 + 24 - 120 + \etc = A
\] 
näherungsweise als
\[
	A = 0 + \frac{1}{e^{9:1}} + \frac{1}{2e^{8:2}} + \frac{1}{3e^{7:3}} + \frac{1}{4e^{6:4}} + \frac{1}{5e^{5:5}} + \frac{1}{6e^{4:6}} + \frac{1}{7e^{3:7}} + \frac{1}{8e^{2:8}} + \frac{1}{9e^{1:9}} + \frac{1}{20} 
\]
erhalten. Diese Terme, weil $e=2,71828128$ ist, werden die folgenden Werte annehmen:
\begin{align*}
	\frac{1}{e^{9:1}} &= 0,00012341 \\
	\frac{1}{2e^{8:2}} &= 0,00915782 \\
	\frac{1}{3e^{7:3}} &= 0,03232399 \\
	\frac{1}{4e^{6:4}} &= 0,05578254 \\
	\frac{1}{5e^{5:5}} &= 0,07357589 \\
	\frac{1}{6e^{4:6}} &= 0,08556952 \\
	\frac{1}{7e^{3:7}} &= 0,09306272 \\
	\frac{1}{8e^{2:8}} &= 0,09735007 \\
	\frac{1}{9e^{1:9}} &= 0,09942659 \\
	\frac{1}{20} &= 0,5000000 \\
	\text{daher}\quad A &= 0,59637255 
\end{align*}
welcher Wert vom wahren schon kaum wahrnehmbar abweicht. Wenn aber die Abszisse in mehrere Teile geteilt worden wäre, dann wäre dieser Wert genauer gefunden worden.
\paragraph{§20}
Weil die Summe als
\[
	A = \int{\frac{e^{1 - 1:x}}{x}\mathrm{d}x}
\]
gefunden worden ist, setze man
\[
	v = e^{1 - 1:x},
\]
sodass für $x=0$ gesetzt auch $v=0$ wird und für $x=1$ gesetzt $v=1$; es wird $1 - \frac{1}{x} = \log{v}$ sein und $x = \frac{1}{1-\log{v}}$ sowie $\log{x} = -\log{(1-\log{v})}$, woher
\[
	\frac{\mathrm{d}x}{x} = \frac{\mathrm{d}v}{v(1-\log{v})}
\]
wird. Weil also
\[
	A = \int{\frac{v\mathrm{d}x}{x}}
\]
ist, nachdem $x=1$ oder $v=1$ gesetzt wurde, wird auch
\[
	A = \int{\frac{\mathrm{d}v}{1 - \log{v}}}
\]
sein, nachdem nach der Integration $v=1$ gesetzt wurde. Es wird aber, indem man die Integration durch eine unendliche Reihe ausführt,
\begin{align*}
	A &= \int{\frac{\mathrm{d}v}{1-\log{v}}} = \frac{v}{1-\log{v}} - \frac{1\cdot v}{(1-\log{v})^2} + \frac{1\cdot 2\cdot v}{(1-\log{v})^3} \\
	&- \frac{1\cdot 2\cdot 3\cdot v}{(1-\log{v})^4} + \frac{1\cdot 2\cdot 3\cdot 4\cdot v}{(1-\log{v})^5} - \etc
\end{align*}
sein und für $v=1$ gesetzt wird wegen $\log{v} = 0$, wie wir angenommen haben,
\[
	A = 1 - 1 + 1\cdot 2 - 1\cdot 2\cdot 3 + 1\cdot 2\cdot 3\cdot 4 - 1\cdot 2\cdot 3\cdot 4\cdot 5 + \etc
\]
sein. Es wird also wiederum $A$ die Fläche der Kurve sein, deren Gestalt zwischen der Abszisse $v$ und der Ordinate $y$ mit dieser Gleichung
\[
	y = \frac{1}{1- \log{v}}
\]
wenn freilich die Abszisse $v=1$ gesetzt wird, in welchem Fall auch $y=1$ wird. Es sollte hier aber bemerkt werden, dass $\log{v}$ den hyperbolischen Logarithmus von $v$ bezeichnet. Nachdem also die Abszisse $v=1$ erneut in $10$ Teile geteilt wurde, werden sich die Ordinaten in den einzelnen Teilungspunkten auf diese Weise verhalten:
\[
\begin{array}{ccl}
\text{wenn $v$ ist} & \qquad & \text{wird $y$ sein} \\
v = \frac{0}{10} & \qquad \qquad & y = 0  \\
v = \frac{1}{10} & \qquad \qquad & y = \frac{1}{(1 + \log{10} - \log{1})} \\
v = \frac{2}{10} & \qquad \qquad & y = \frac{1}{(1 + \log{10} - \log{2})} \\
v = \frac{3}{10} & \qquad \qquad & y = \frac{1}{(1 + \log{10} - \log{3})} \\
v = \frac{4}{10} & \qquad \qquad & y = \frac{1}{(1 + \log{10} - \log{4})} \\
v = \frac{5}{10} & \qquad \qquad & y = \frac{1}{(1 + \log{10} - \log{5})} \\
v = \frac{6}{10} & \qquad \qquad & y = \frac{1}{(1 + \log{10} - \log{6})} \\
v = \frac{7}{10} & \qquad \qquad & y = \frac{1}{(1 + \log{10} - \log{7})} \\
v = \frac{8}{10} & \qquad \qquad & y = \frac{1}{(1 + \log{10} - \log{8})} \\
v = \frac{9}{10} & \qquad \qquad & y = \frac{1}{(1 + \log{10} - \log{9})} \\
v = \frac{10}{10} & \qquad \qquad & y = \frac{1}{1}
\end{array}
\]
Und daher wird man durch Annäherung der Fläche den Wert des Buchstaben $A$ hinreichend genau erhalten.
\paragraph{§21}
Es ist aber eine andere Art die Summe dieser Reihe zu untersuchen gegeben, die man aus der Lehre der Kettenbrüche herholt und die um Vieles leichter und schneller die Aufgabe erledigt; es sei nämlich, indem man die Formel allgemeiner ausdrückt
\[
	A = 1 - 1x + 2x^2 - 6x^3 + 24x^4 - 120x^5 + 720x^6 - 5040x^7 + \etc = \frac{1}{1+B};
\] 
es wird
\[
	B = \frac{1x - 2x^2 + 6x^3 - 24x^4 + 120x^5 - 720x^6 + 5040x^7 - \etc}{1 - 1x + 2x^2 - 6x^3 + 24x^4 - 120x^5 + 720x^6 - 5040x^7 + \etc} = \frac{x}{1+C}
\]
sein und
\[
	1 + C = \frac{1 - 1x + 2x^2 - 6x^3 + 24x^4 - 120x^5 + 720x^6 - 5040x^7 + \etc}{1 - 2x + 6x^2 - 24x^3 + 120x^4 - 720x^5 + 5040x^6 - \etc}.
\]
Also
\[
	C = \frac{x - 4x^2 + 18x^3 - 96x^4 + 600x^5 - 4320x^6 + \etc}{1 - 2x + 6x^2 - 24x^3 + 120x^4 - 720x^5 + \etc} = \frac{x}{1+D}
\]
daher
\[
	D = \frac{2x - 12x^2 + 72x^3 - 480x^4 + 3600x^5 - \etc}{1 - 4x + 18x^2 - 96x^3 + 600x^4 - \etc} = \frac{2x}{1+E}
\]
Weiter
\[
	E = \frac{2x - 18x^2 + 144x^3 - 1200x^4 + \etc}{1 - 6x + 36x^2 - 240x^3 + \etc} = \frac{2x}{1-F}
\]
und
\[
	F = \frac{3x - 36x^2 + 360x^3 - \etc}{1 - 9x + 72x^2 - 600x^3 + \etc} = \frac{3x}{1+G}
\]
Es wird
\[
	G = \frac{3x - 48x^2 + \etc}{1 - 12x + 120x^2 - \etc}  =\frac{3x}{1+H}
\]
sein. So
\[
	H = \frac{4x - \etc}{1 - 16x + \etc} = \frac{4x}{1 + I}
\]
und so weiter wird klar sein, dass
\[
	I = \frac{4x}{1 + K},\quad K = \frac{5x}{1+L},\quad L = \frac{5x}{1+M},\quad \text{etc ins Unendliche}
\]
sein wird, sodass man die Struktur dieser Formeln leicht durchschaut. Nachdem aber diese Werte nacheinander eingesetzt wurden, wird
\[
	1 - 1x + 2x^2 - 6x^3 + 24x^4 - 120x^5 + 720x^6 - 5040x^7 + \etc
\]
gleich
\[
	A = \cfrac{1}{1 + \cfrac{x}{1 + \cfrac{x}{1 + \cfrac{2x}{1 + \cfrac{2x}{1 + \cfrac{3x}{1 + \cfrac{3x}{1 + \cfrac{4x}{1 + \cfrac{4x}{1 + \cfrac{5x}{1 + \cfrac{5x}{1 + \cfrac{6x}{1 + \cfrac{6x}{1 + \cfrac{7x}{\etc.}}}}}}}}}}}}}}
\]
sein.
\paragraph{§22}
Wie aber der Wert von Kettenbrüchen solcher Art zu untersuchen ist, habe ich anderenorts gezeigt. Weil natürlich die ganzen Teile der einzelnen Nenner Einheiten sind, gehen allein die Zähler in die Rechnung ein; es sei also $x=1$ und die Untersuchung der Summe $A$ wird auf folgende Weise angestellt werden:
\[
\begin{array}{rclllllllllll}
A &= &\frac{0}{1}, &\frac{1}{1}, &\frac{1}{2}, &\frac{2}{3}, &\frac{4}{7}, &\frac{8}{13}, &\frac{20}{34}, &\frac{44}{73}, &\frac{124}{209}, &\frac{300}{501}, &\etc \\[1ex]
\text{Zähler}:& &1, &1, &2, &2, &3, &3, &4, &4, &5, &5, &\etc
\end{array}
\]
Die hier beschaffenen Brüche nähern sich natürlich immer besser dem wahren Wert von $A$ und sind freilich größer und kleiner als dieser, sodass
\begin{alignat*}{6}
	A > \frac{0}{1},\quad & A > \frac{1}{2},\quad & A > \frac{4}{7},\quad & A > \frac{20}{34},\quad & A > \frac{124}{209},\quad & \etc \\
	A < \frac{1}{1},\quad & A < \frac{2}{3},\quad & A < \frac{8}{13},\quad & A < \frac{44}{73},\quad & A < \frac{300}{501},\quad & \etc
\end{alignat*}
ist. Daher werden in Dezimalbrüchen die Werte von $A$
\[
\begin{array}{ccc}
\text{zu kleine Werte} & \qquad & \text{zu große Werte} \\
0,0000000000 & \qquad \qquad & 1,0000000000 \\
0,5000000000 & \qquad \qquad & 0,6666666667 \\
0,5714285714 & \qquad \qquad & 0,6153846154 \\
0,5882352941 & \qquad \qquad & 0,6027397260 \\
0,5933001436 & \qquad \qquad & 0,5988023952 
\end{array}
\]
sein. Wenn gleich zwischen den zu großen und zu kleinen Näherungstermen die arithmetischen Mittel genommen werden, werden erneut abwechselnd zu große und zu kleine Werte hervorgehen, welche die folgenden sein werden:
\[
\begin{array}{ccc}
\text{zu kleine Werte} & \qquad & \text{zu große Werte} \\
0,5000000000 & \qquad \qquad & 0,7500000000 \\
0,5833333333 & \qquad \qquad & 0,6190476190 \\
0,5934065934 & \qquad \qquad & 0,6018099548 \\
0,5954875100 & \qquad \qquad & 0,5980205807 \\
0,5960519153 & \qquad \qquad &  
\end{array}
\]
und so gelangen wir schon hinreichend nahe an den wahren Wert von $A$.
\paragraph{§23}
Wir können aber den Wert dieses unendlichen Bruches durch Teile auf die Art untersuchen: Es sei
\[
	A = \cfrac{1}{1 + \cfrac{1}{1+ \cfrac{1}{1 + \cfrac{2}{1 + \cfrac{2}{1 + \cfrac{3}{1 + \cfrac{3}{1 + \cfrac{4}{1 + \cfrac{4}{1 + \cfrac{5}{1 + \cfrac{5}{1 + \cfrac{6}{1 + \cfrac{6}{1 + \cfrac{7}{1 + \cfrac{7}{1 + \cfrac{8}{1 + \cfrac{8}{1 + p}}}}}}}}}}}}}}}}}
\]
und
\[
	p = \cfrac{9}{1 + \cfrac{9}{1 + \cfrac{10}{1 + \cfrac{10}{1 + \cfrac{11}{1 + \cfrac{11}{1 + \cfrac{12}{1 + \cfrac{12}{1 + \cfrac{13}{1+ \cfrac{13}{1 + \cfrac{14}{1 + \cfrac{14}{1 + \cfrac{15}{1 + \cfrac{15}{1 +q}}}}}}}}}}}}}}
\]
und
\[
	q = \cfrac{16}{1 + \cfrac{16}{1 + \cfrac{17}{1 + \cfrac{17}{1 + \cfrac{18}{1+ \cfrac{18}{1 + \cfrac{19}{1 + \cfrac{19}{1 + \cfrac{20}{1 + \cfrac{20}{1 + r}}}}}}}}}}
\]
es wird
\[
	r = \cfrac{21}{1 + \cfrac{21}{1 + \cfrac{22}{1 + \cfrac{22}{1 + \cfrac{23}{1 + \cfrac{23}{1 + \etc}}}}}}
\]
Nachdem diese Werte entwickelt wurden, wird man zuerst
\[
	A = \frac{491459820 + 139931620p}{824073141 + 234662231p}
\]
finden, darauf
\[
	p = \frac{2381951 + 649286q}{887640 + 187440q}
\]
und
\[
	q = \frac{11437136 + 2924816r}{3697925 + 643025r}.
\]
Es ist also übrig, dass der Wert von $r$ bestimmt wird, was freilich genauso schwer wie der von $A$ ist, aber es genügt hier, den Wert von $r$ nur näherungsweise zu kennen; der Fehler nämlich, der beim Wert von $r$ begangen wurde, bewirkt einen um Vieles kleineren Fehler beim Wert von $q$ und daher besiedelt ein erneut weit kleinerer Fehler den Wert von $p$; daraus schließlich wird der Fehler, der den Wert von $A$ befleckt im Ganzen unentdeckt bleiben.
\paragraph{§24}
Weil darauf die Zähler $21$, $21$, $22$, $22$, $23$, $23$, etc, die in den Kettenbruch $r$ eingehen, schon nähernd an das Verhältnis der Gleichheit herangehen, zumindest anfangs, kann daher Hilfe geholt werden um einen Wert besser zu erkennen. Wenn nämlich diese Zähler alle gleich wären, dass
\[
	r = \cfrac{21}{1 + \cfrac{21}{1 + \cfrac{21}{1 + \cfrac{21}{1 + \etc}}}}
\]
wäre, würde
\[
	r = \frac{21}{1+r}
\]
sein und daher
\[
	rr+r = 21
\]
und 
\[
	r = \frac{\sqrt{85}-1}{2}.
\]
Weil aber diese Nenner wachsen, wird dieser Wert größer sein als der rechtmäßige. Dennoch lässt sich schließen, wenn die drei folgenden Kettenbrüche
\begin{minipage}{0.5\textwidth}
\centering
\[
	r = \cfrac{21}{1 + \cfrac{21}{1 + \cfrac{22}{1 + \cfrac{22}{1 + \cfrac{23}{1 + \etc}}}}}
\]
\end{minipage}
\begin{minipage}{0.5\textwidth}
\centering
\[
	s = \cfrac{22}{1 + \cfrac{22}{1 + \cfrac{23}{1 + \cfrac{23}{1 + \cfrac{24}{1 + \etc}}}}}
\]
\end{minipage}
\begin{center}
\[
	t = \cfrac{23}{1 + \cfrac{23}{1 + \cfrac{24}{1 + \cfrac{24}{1 + \cfrac{25}{1 + \etc}}}}}
\]
\end{center}
festgesetzt werden, dass die Werte der Größen $r$, $s$, $t$ in einer arithmetischen Progression fortschreiten und dass $r+s = 2t$ sein wird; daher wird der Wert von $r$ hinreichend genau berechnet werden; Damit sich aber diese Untersuchung weiter erstreckt, wollen wir für die Zahlen $21$, $22$, $23$ diese unbestimmten $a-1$, $a$ und $a+1$ annehmen, dass

\begin{minipage}{0.5\textwidth}
\centering
\[
	r = \cfrac{a-1}{1 + \cfrac{a-1}{1 + \cfrac{a}{1 + \cfrac{a}{1 + \cfrac{a+1}{1 + \etc}}}}}
\]
\end{minipage}
\begin{minipage}{0.5\textwidth}
\centering
\[
	s = \cfrac{a}{1 + \cfrac{a}{1 + \cfrac{a+1}{1 + \cfrac{a+1}{1 + \cfrac{a+2}{1 + \etc}}}}}
\]
\end{minipage}
\begin{center}
\[
	t = \cfrac{a+1}{1 + \cfrac{a+1}{1 + \cfrac{a+2}{1 + \cfrac{a+2}{1 + \cfrac{a+3}{1 + \etc}}}}}
\]
\end{center}
ist und es wird
\[
	r = \cfrac{a-1}{1 + \cfrac{a-1}{1+s}} \qquad s = \cfrac{a}{1 + \cfrac{a}{1+t}}
\]
sein, woher
\[
	r = \frac{(a-1)s + a -1 }{s+a}
\]
bewirkt wird und
\[
	s = \frac{at+a}{t+a+1} \quad \text{oder} \quad t = \frac{(a+1)s-a}{a-s},
\]
woher
\[
	r+t = \frac{2ss + (2aa - 2a + 1)s - a}{aa - ss} = 2s
\]
wird; und daher wird
\[
	2s^3 + 2ss - (2a-1)s - a = 0
\]
sein, aus welcher Gleichung sich der Wert von $s$ und daher weiter der Wert von $r$ bestimmen lässt.
\paragraph{§25}
Es sei nun $a=22$ und wir werden diese zu lösende kubische Gleichung haben
\[
	2s^3 + 2ss - 43s - 22 = 0,
\]
deren Wurzel sofort zwischen den Grenzen $4$ und $5$ liegend entdeckt wird. Es sei daher $s=4+u$ und es wird
\[
	34 = 69u + 26uu + 2u^3
\]
sein. Es sei weiter $u=0,4 + v$; es wird
\[
	u^2 = 0,16 + 0,8v + vv \quad \text{und} \quad u^3 = 0,064 + 0,48v + 1,2v^2 + v^3
\]
sein und daher
\[
	2,112 = 90,76v + 28,4v^2 + 2v^3,
\]
woher näherungsweise
\[
	v = 0,023 \quad \text{und} \quad s = 4,423
\]
sein wird. Weil also
\[
	r = \frac{21s+21}{s+2}
\]
ist, wird
\[
	r = \frac{113,883}{26,423} = 4,31
\]
sein und daher weiter
\[
	q = \frac{24043093}{6469363} = 3,71645446,
\]
woher man
\[
	p = \frac{4794992,85}{1584252,22} = 3,0266600163
\]
erhält und daher schließlich
\[
	A = \frac{914985259,27}{1534315932,90} = 0,5963473621372,
\]
welcher Wert in einen Kettenbruch verwandelt
\[
	A = \cfrac{1}{1 + \cfrac{1}{1+ \cfrac{1}{2 + \cfrac{1}{10 + \cfrac{1}{1 + \cfrac{1}{1 + \cfrac{1}{4 + \cfrac{1}{2 + \cfrac{1}{2 + \cfrac{1}{13 + \cfrac{1}{4 + \etc}}}}}}}}}}}
\]
gibt, woher man die folgenden Brüche findet, die den Wert von $A$ näherungsweise beschaffen
\[
\begin{array}{cccccccccccl}
	& 1 & 1 & 2 & 10 & 1 & 1 & 4 & 2 & 2 & 13\\[1ex]
	A = & \frac{0}{1}, &\frac{1}{1}, & \frac{1}{2}, & \frac{3}{5}, & \frac{31}{52}, & \frac{34}{57}, & \frac{65}{109}, & \frac{294}{493}, & \frac{653}{1095}, & \frac{1600}{2683}, & \etc 
\end{array}
\]
Diese Brüche aber sind abwechselnd größer und kleiner als der Wert von $A$ und der letzte freilich $\frac{1600}{2683}$ ist zu groß, die Abweichung ist dennoch kleiner als $\frac{1}{2683 \cdot 35974}$; weil daher
\[
	\frac{1}{A} = \frac{2683}{1600}
\]
ist, wird näherungsweise
\[
	\frac{1}{A} = 1,676875
\]
sein.
\paragraph{§26}
Die Methode, die ich oben in §$21$ benutzt habe, um diese Reihe
\[
	1 - 1x + 2x^2 - 6x^3 + 24x^4 - 120x^5 + 720x^6 - 5040x^7 + \etc
\]
in einen Kettenbruch zu verwandeln, erstreckt sich weiter und kann auf ähnliche Art auf diese um Vieles angenehmere Reihe angewandt weden:
\begin{align*}
	z &= 1 -mx +m(m+n)x^2 - m(m+n)(m+2n)x^3 \\
	&+m(m+n)(m+2n)(m+3n)x^4 - \etc;
\end{align*}
man wird nämlich, nachdem dieselben Operationen angestellt wurden,
\[
	z = \cfrac{1}{1 + \cfrac{mx}{1 + \cfrac{nx}{1 + \cfrac{(m+n)x}{1 + \cfrac{2nx}{1 + \cfrac{(m+2n)x}{1 + \cfrac{3nx}{1 + \cfrac{(m+3n)x}{1 + \cfrac{4nx}{1 + \cfrac{(m+4n)x}{1 + \cfrac{5nx}{1 + \etc}}}}}}}}}}}
\]
finden. Derselbe Ausdruck aber und andere ähnliche können leicht mithilfe der Theoreme gefunden werden, die ich in meinen Abhandlungen über Kettenbrüche in \textit{Comment. Acad. Petropol.} gezeigt habe. Ich habe nämlich gezeigt, dass dieser Gleichung
\[
	ax^{m-1}\mathrm{d}x = \mathrm{d}z + cx^{n-m-1}z\mathrm{d}x + bx^{n-1}z\mathrm{d}x
\]
dieser Wert
\[
	z = \cfrac{ax^m}{m + \cfrac{(ac+mb)x^n}{m + n + \cfrac{(ac-nb)x^n}{m + 2n + \cfrac{(ac + (m+n)b)x^n}{m + 3n + \cfrac{(ac-2nb)x^n}{m + 4n + \cfrac{(ac + (m+2n)b)x^n}{m + 5n + \cfrac{(ac-3nb)x^n}{m + 6n + \etc}}}}}}}
\]
genügt. Wenn also $c=0$ ist, wird
\[
	\mathrm{d}z + bx^{n-1}z\mathrm{d}x = ax^{m-1}\mathrm{d}x
\]
sein und
\[
	e^{bx^n : n}z = a\int{e^{bx^n : n}x^{m-1}\mathrm{d}x} \quad \text{und} \quad z = ae^{-bx^n : n}\int{e^{bx^n : n}x^{m-1}\mathrm{d}x}
\]
und durch eine Reihe
\[
	z = \frac{ax^m}{m} - \frac{abx^{m+n}}{m(m+n)} + \frac{ab^2 x^{m+2n}}{m(m+n)(m+2n)} - \frac{ab^3 x^{m+3n}}{m(m+n)(m+2n)(m+3n)} + \etc.
\]
In dieser aber ist unsere Form, die wir betrachten, nicht enthalten.
\paragraph{§27}
Ich habe aber weiter gefunden, wenn man diese Gleichung hat
\[
	fx^{m+n}\mathrm{d}x = x^{m+1}\mathrm{d}z + ax^{m}z\mathrm{d}x + bx^n z\mathrm{d}x + czz\mathrm{d}x,
\]
dass der Wert von $z$ durch einen unendlichen Bruch dieser Art ausgedrückt wird
\[
	z = \cfrac{fx^m}{b + \cfrac{(mb+ab+cf)x^{m-n}}{b + \cfrac{(mb - nb + cf)x^{m-n}}{b + \cfrac{(2mb-nb+ab+cf)x^{m-n}}{b + \cfrac{(2mb-2nb+cf)x^{m-n}}{b + \cfrac{(3mb-2nb+ab+cf)x^{m-n}}{b + \cfrac{(3mb-3nb+cf)x^{m-n}}{b + \etc}}}}}}}
\]
Damit wir also denselben Wert $z$ angenehm durch eine gewöhnliche Reihe ausdrücken können, sei $c=0$, dass man diese Gleichung hat
\[
	fx^{m+n}\mathrm{d}x = x^{m+1}\mathrm{d}z + ax^m z\mathrm{d}x + bx^n z\mathrm{d}x,
\]
und es wird durch einen Kettenbruch
\[
	z = \cfrac{fx^m}{b + \cfrac{b(m+a)x^{m-n}}{b + \cfrac{b(m-n)x^{m-n}}{b + \cfrac{b(2m-n+a)x^{m-n}}{b + \cfrac{b(2m-2n)x^{m-n}}{b + \cfrac{b(3m-2n+a)x^{m-n}}{b + \cfrac{b(3m-3n)x^{m-n}}{b + \etc}}}}}}}
\]
sein. Durch Integrieren aber wird
\[
	x^a e^{bx^{n-m}:(n-m)}z = f\int{e^{bx^{n-m}:(m-n)}}x^{a+n-1}\mathrm{d}x
\]
sein oder, wenn $m-n=k$ ist, wird
\[
	z = fe^{b:kx^k}x^{-a}\int{e^{-b:kx^k} x^{a+n-1}\mathrm{d}x}
\]
sein, wenn freilich die Integration so ausgeführt wird, dass $z$ für $x=0$ gesetzt verschwindet. Durch eine unendliche Reihe wird aber
\begin{align*}
	z =& \frac{f}{b}x^m - \frac{(m+a)}{b^2}fx^{2m-n} + \frac{(m+a)(2m-n+a)f}{b^3} x^{3m-2n} \\
	&-\frac{(m+a)(2m-n+a)(3m-2n+a)f}{b^4}x^{4m-3n} \\
	&+\frac{(m+a)(2m-n+a)(3m-2n+a)(4m-3n+a)f}{b^5}x^{5m-4n} - \etc
\end{align*}
sein.
\paragraph{§28}
Damit diese Ausdrücke einfach werden und trotzdem nicht deren Verallgemeinerung eingeschränkt wird, setze man
\[
	b = 1, \quad f = 1, \quad m+a=p, \quad m-n = q,
\]
dass
\[
	a = p-m \quad \text{und} \quad n = m-q
\]
ist; man wird diese Differentialgleichung haben
\[
	x^m\mathrm{d}x = x^{q+1}\mathrm{d}z + (p-m)x^q z\mathrm{d}x + z\mathrm{d}x,
\]
deren Integral zuerst
\[
	z = e^{1:qx^q}x^{m-p}\int{e^{-18qx^q}x^{p-q-1}\mathrm{d}x}
\]
ist. Derselbe Wert der Größe $z$ wird weiter durch die folgende unendliche Reihe ausgedrückt werden
\[
	z = x^m - px^{m+q} + p(p+q)x^{m+2q} - p(p+q)(p+2q)x^{m+3q} + \etc.
\]
Schließlich wird dieser Reihe dieser Kettenbruch
\[
	z = \cfrac{x^m}{1 + \cfrac{px^q}{1 + \cfrac{qx^q}{1 + \cfrac{(p+q)x^q}{1 + \cfrac{2qx^q}{1 + \cfrac{(p+2q)x^q}{1 + \cfrac{3qx^q}{1 + \cfrac{(p+3q)x^q}{1 + \etc}}}}}}}}
\]
äquivalent sein, welcher Ausdruck natürlich mit dem, den wir zuvor in §$26$ erhalten haben, übereinstimmt; und weil ja über die Art und Weise, auf die wir jenen gefunden haben, noch gezweifelt werden könnte, ob die Zähler nach dem beobachteten Gesetz ins Unendliche fortschreiten oder nicht, ist dieser Zweifel nun völlig ausgeräumt worden. Es liefert nämlich diese Betrachtung eine sichere Methode, unzählige divergente Reihen zu summieren oder selbigen äquivalente Werte zu finden; unter diesen ist der, den wir betrachtet haben, ein Spezialfall.
\paragraph{§29}
Es scheint aber weiter der Fall bemerkenswert, in dem $p=1$ und $q=2$ sowie $m=1$ ist; es wird nämlich
\[
	z = e^{1:2xx}\int{e^{-1:2xx}}\mathrm{d}x : xx
\]
sein und die unendliche Reihe wird sich so verhalten
\[
	z = x - 1x^3 + 1\cdot 3x^5 - 1\cdot 3\cdot 5x^7 + 1\cdot 3\cdot 5\cdot 7x^9 - \etc,
\]
die diesem Kettenbruch gleich ist
\[
	z = \cfrac{x}{1 + \cfrac{1xx}{1 + \cfrac{2xx}{1 + \cfrac{3xx}{1 + \cfrac{4xx}{1 + \cfrac{5xx}{1 + \cfrac{6xx}{1 + \etc}}}}}}}
\]
Wenn deshalb $x=1$ gesetzt wird, dass 
\[
	z = 1 - 1 + 1\cdot 3 - 1\cdot 3\cdot 5 + 1\cdot 3\cdot 5\cdot 7 - 1\cdot 3\cdot 5\cdot 7\cdot 9 + \etc
\]
wird, welche Reihe besonders divergent ist, kann ihr Wert trotzdem durch diesen konvergenten Kettenbruch ausgedrückt werden
\[
	z = \cfrac{1}{1 + \cfrac{1}{1 + \cfrac{2}{1 + \cfrac{3}{1 + \cfrac{4}{1 + \cfrac{5}{1 + \etc}}}}}}
\]
der die folgenden dem wahren Wert von $z$ näherungsweise gleiche Brüche liefert
\[
\begin{array}{cccccccccccccc}
	& 1 & 2 & 3 & 4 & 5 & 6 & 7 & 8 & 9 & 10 & 11 & 12& \\[1ex]
	z = & \frac{0}{1}, & \frac{1}{1}, & \frac{1}{2}, & \frac{3}{4}, & \frac{6}{10}, & \frac{18}{26}, & \frac{48}{76}, & \frac{156}{232}, &\frac{492}{764}, & \frac{1740}{2620}, & \frac{6168}{9496}, & \frac{23568}{35696} & \etc;
\end{array}
\]
wenn also
\[
	z = \cfrac{1}{1 + \cfrac{1}{1 + \cfrac{2}{1 + \cfrac{3}{1 + \cfrac{4}{1 + \cfrac{5}{1 + \cfrac{6}{1 + \cfrac{7}{1 + \cfrac{8}{1 + \cfrac{9}{1 + \cfrac{10}{1 + p}}}}}}}}}}}
\]
ist, wird
\[
	z = \frac{23568 + 6168p}{35696 + 9496p}
\]
sein oder
\[
	z = \frac{2946 + 771p}{4402 + 1187p}
\]
und
\[
	p = \cfrac{11}{1 + \cfrac{12}{1 + \cfrac{13}{1 + \cfrac{14}{1 + \cfrac{15}{1 + \etc}}}}}
\]
Es sei
\[
	p = \frac{11}{1+q} \quad \text{und} \quad q = \frac{12}{1 + r}
\]
und weil $p$, $q$, $r$ gleichmäßig wachsen, wird
\[
	2q = \frac{12 + 22q -qq}{q+qq} \quad \text{und} \quad 2q^3 + 3qq - 22q - 12 = 0
\]
sein, wo näherungsweise
\[
	q = 2,94,\quad p = 2,79 \quad \text{und} \quad z = \frac{5097,09}{7773,73} = 0,65568
\]
ist.
\end{document}